\documentclass{elsart}
\usepackage{amssymb}
\usepackage{amsfonts}
\newcommand{\trace}{\mathop{\rm Tr}\nolimits}

\newcommand{\twomat}[4]{\left(\begin{array}{cc}#1&#2\\#3&#4\end{array}\right)}

\newcommand{\cS}{{\mathcal S}}

\newcommand{\R}{{\mathbb{R}}}

\DeclareRobustCommand\openone{\leavevmode\hbox{\small1\normalsize\kern-.33em1}}
\newcommand{\id}{\mathrm{\openone}}

\newcommand{\be}{\begin{equation}}
\newcommand{\ee}{\end{equation}}
\newcommand{\bea}{\begin{eqnarray}}
\newcommand{\eea}{\end{eqnarray}}
\newcommand{\beas}{\begin{eqnarray*}}
\newcommand{\eeas}{\end{eqnarray*}}
\newtheorem{definition}{Definition}
\newtheorem{theorem}{Theorem}
\newtheorem{lemma}{Lemma}

\begin{document}
\begin{frontmatter}
\title{In-betweenness, a Geometrical Monotonicity Property for Operator Means}
\author{Koenraad M.R.\ Audenaert}
\address{
Dept.\ of Mathematics,\\
Royal Holloway, University of London,\\
Egham TW20 0EX, United Kingdom}
\ead{koenraad.audenaert@rhul.ac.uk}
\date{\today}
\begin{abstract}
We introduce the notions of in-betweenness and monotonicity with respect to a metric for operator means.
These notions can be seen as generalising their natural counterpart for scalar means, and
as a relaxation of the notion of geodesity. We exhibit two classes of non-trivial means that are monotonic
with respect to the Euclidean metric. We also show that all Kubo-Ando means are monotonic with respect to
the trace metric, which is the natural metric for the geometric mean.
\end{abstract}

\end{frontmatter}
\section{Introduction}
According to the highly respected Merriam-Webster's dictionary, a \textit{mean}
is ``\textit{a value that lies within a range of values and is computed according to a prescribed law.}''
The best-known examples of means in this sense are the arithmetic mean and the geometric mean of two real scalars
$x$ and $y$, defined by the `prescribed laws' $\mu(x,y)=(x+y)/2$ and $\mu(x,y)=\sqrt{xy}$, respectively.
As is easily checked, these means indeed lie `within the range' $[x,y]$. 
Many more means have been defined, like the harmonic mean
$\mu(x,y)=2(1/x+1/y)^{-1}$
and the power means $\mu(x,y)=((x^p+y^p)/2)^{1/p}$ (with $p\ge1$), 
and they all share this property of being contained
in the interval $[x,y]$. At least for real numbers, the dictionary definition appears mathematically correct.
For succinctness, we will call this property that for all $x\le y$, $x\le \mu(x,y)\le y$, 
the \textit{in-betweenness} property of a mean $\mu$. 

The basic notion of mean has been extended to more general mathematical objects,
like functions, vectors, matrices and operators. 
Because of the more complicated structure of these objects,
it no longer makes
sense in general to say that the mean of objects $f$ and $g$ `lies within a range' defined by $f$ and $g$.
The definition of in-betweenness for scalar means inherently relies on the endowment of $\R$ with
a total ordering, $(\R,\le)$. For more complicated structures a partial ordering is the best one can hope for,
which in itself does not provide a solid foundation for an in-betweenness property. 

In a number of cases the geometry of the space in which the mean is defined induces a total ordering;
this happens when the mean can be parameterised as $t\mapsto \mu(X,Y,t)$ and the path traced out by
varying $t$ is a geodesic with respect to the chosen metric of the space. 
A well-known example is the geometric mean,
which can be parameterised as $x\#_t y = x^t y^{1-t}$, or as $A\#_t B=A^{1/2}(A^{-1/2}BA^{1/2})^t A^{1/2}$ for
positive operators. It can be shown that the path $t\mapsto A\#_t B$ is a geodesic with
respect to the trace metric \cite{bhatia_EMI} (see below). 
In-betweenness with respect to the metric then follows by definition.

In general, however, it need not be straightforward to parameterise a given mean and then find a metric 
such that the mean lies on a geodesic. Secondly, the context in which the mean is to be used might impose
a different metric and checking in-betweenness is no longer trivial.
Thirdly, it is fair to say that most means have not been defined starting from such geometric considerations.
Often, the only claim that such means can lay on their name is the close resemblance between
their defining prescribed law and a similar law defined for their scalar counterpart.
The fact that one has proceeded with the definition of these means anyway is largely due to their applicability.
Amongst the better-known means for positive operators 
are the \textit{arithmetic mean} $(A,B)\mapsto (A+B)/2$,
the \textit{geometric mean} $(A,B)\mapsto A\# B = A^{1/2}\,\,(A^{-1/2} B A^{-1/2})^{1/2}\,\,A^{1/2}$
and the \textit{harmonic mean} $(A,B)\mapsto A!B = ((A^{-1}+B^{-1})/2)^{-1}$ (for invertible $A$ and $B$, that is).

Within the last three decades, the area of operator means became largely dominated
by what is now known as the class of Kubo-Ando operator means. 
In a beautiful and very influential paper,
Kubo and Ando \cite{kuboando} introduced a set of axioms and showed that they were satisfied by a large number of the
then known operator means, including the abovementioned arithmetic, geometric and harmonic mean.
Moreover, they completely characterised the class of such means and showed that they
are in one-to-one correspondence with the non-negative operator monotone functions on $(0,+\infty)$.
The ensuing theory stimulated a lot of research because 
of its connections to Riemannian geometry, and its applications in mathematical physics.
As a case in point, one should note that the Kubo-Ando axioms do not appeal to any underlying
geometry, metric or geodesic.

It has to be emphasised that not all operator means in current use are Kubo-Ando means.
We mention only two prominent examples here, as they are the subject of the technical part of this paper.
Our first example is the class of \textit{power means}, studied for operators 
by Bhagwat and Subramanyan \cite{bhagwat}.
They are defined as
\be\label{eq:power}
(A,B) \mapsto ((A^p+B^p)/2)^{1/p},
\ee
with $p\in \R$.
Clearly, this class contains the arithmetic mean ($p=1$), and the harmonic mean ($p=-1$).
These power means are Kubo-Ando means only when $-1\le p\le 1$.
In spite of this, the power means with $p>1$ have many important applications, e.g.\ in mathematical physics 
and in the theory of operator spaces, where they form the basis of certain generalisations of
$\ell_p$ norms to non-commutative vector-valued $L_p$ spaces \cite{CL99}.

Our second example of non-Kubo-Ando means is the class of Heinz means.
The Heinz means for non-negative scalars are 
weighted versions of the geometric mean:
$H_\nu(x,y) = x^\nu y^{1-\nu}$, with $0\le\nu\le 1$.
Sometimes another definition is adopted that is slightly more symmetrical \cite{bhatia_heinz}.
Namely, the symmetric Heinz mean is defined as 
$H'_\nu(x,y) = (x^\nu y^{1-\nu} + x^{1-\nu} y^\nu)/2$, which is invariant under replacing $\nu$ by $1-\nu$.
The reason for this convention is that the symmetric Heinz mean interpolates between the arithmetic mean
($H'_0(x,y)=H'_1(x,y) = (x+y)/2$) and the geometric mean ($H'_{1/2}(x,y)=\sqrt{xy}$).

These two definitions carry over to operators in a straightforward way: one defines the symmetric Heinz mean
as $H'_\nu(A,B) = (A^\nu B^{1-\nu} + A^{1-\nu} B^\nu)/2$, and the unsymmetric one as $H_\nu(A,B) = A^\nu B^{1-\nu}$.
Clearly, these means cannot be Kubo-Ando means as they violate the first axiom of closure. In general,
the Heinz mean of two positive operators is not even self-adjoint, let alone positive.

Nevertheless, the Heinz means have great importance.
The unsymmetric Heinz mean, in particular, is a basic quantity in quantum physics. 
When applied to density operators, the logarithm of the
trace $\log\trace H_\nu(\rho,\sigma) = \log\trace \rho^\nu \sigma^{1-\nu}$ is known as the relative Renyi entropy.
The normalised mean itself, $\rho^\nu\sigma^{1-\nu}$ divided by its trace,
can be considered as a quantum generalisation of the so-called Hellinger arc between two probability distributions
\cite{chernoff}.
In the present manuscript we will only consider the unsymmetric version of the Heinz mean, 
for reasons of simplicity.

In this paper we shall investigate one possible route towards defining an in-betweenness property
for operator means, overcoming the lack of a total ordering and of the existence of a natural metric. 
In fact we will define two varieties of such a property. Both are based on endowing the set of positive
operators with a simple Euclidean geometry; this is the topic of Section \ref{sec:euclidean}.
For any definition to be meaningful, one would normally expect the existence of objects that satisfy it.
We show that there exist indeed non-trivial operator means (apart from the arithmetic mean) 
that satisfy this kind of in-betweenness,
namely the power means, in Section \ref{sec:power}, and the Heinz means, in Section \ref{sec:heinz}. 
Next, in Section \ref{sec:KA} 
we will exhibit a counterexample that shows that Kubo-Ando means generally do not satisfy
in-betweenness with respect to the Euclidean distance.
In contrast, we will prove that they all satisfy in-betweenness with respect to the trace metric distance,
the metric whose geodesics are traced out by the geometric means.
In Section \ref{sec:discussion} we conclude and briefly state further research directions.
\section{Kubo-Ando means}
The Kubo-Ando axioms are the following,
with $\sigma$ the generic symbol for a mean in the Kubo-Ando sense,
and $A,B,C,D$ arbitrary non-negative operators:
\begin{enumerate}
\item Closure: A mean is a binary operation on and into the class of positive operators, $A\sigma B\ge0$;
\item Monotonicity: $A\le C$ and $B\le D$ imply $A\sigma B\le C\sigma D$;
\item Transformer inequality: $C(A\sigma B)C \le (CAC)\sigma(CBC)$;
\item Continuity: $A_n\downarrow A$ and $B_n\downarrow B$ imply $(A_n\sigma B_n)\downarrow (A\sigma B)$;
\item Normalisation: $\id\sigma\id=\id$.
\end{enumerate}
Here, the notation $A_n\downarrow A$ is a shorthand for the statement that there is a sequence of positive operators
$A_1\ge A_2\ge\ldots A_n$ with $A_n$ converging strongly to $A$. For further information about these
axioms we refer to \cite{kuboando}.

Dropping the normalisation condition, Kubo and Ando then showed that for any mean $\sigma$
satisfying these axioms, the function $f(x):=1\sigma x$ is a non-negative operator monotone
function on $(0,+\infty)$. Conversely, for any
non-negative operator monotone function $f$ on $(0,+\infty)$,
there is a mean satisfying the axioms, via the construction
$$
A\sigma B = A^{1/2}\,\, f(A^{-1/2} B A^{-1/2})\,\, A^{1/2}.
$$
Because of this correspondence, $f$ is called the representing function of the mean.

For example, the power means (\ref{eq:power}) 
are Kubo-Ando means only when $-1\le p\le 1$, as this is the condition for operator monotonicity
of the representing function $f(x)=((1+x^p)/2)^{1/p}$.


Exploiting the theory of operator monotone functions, Kubo and Ando arrived at an integral representation
of any mean satisfying their axioms (excluding the normalisation condition), see Theorem 3.4 in \cite{kuboando}.
Given any positive Radon measure $\mu(s)$ on $[0,+\infty]$, there is a unique corresponding Kubo-Ando mean
represented as
\be\label{eq:KAint1}
A\sigma B = aA+bB+\int_{(0,+\infty)} \frac{1+s}{s}\,\,(sA):B\,\,d\mu(s),
\ee
with $a=\mu(\{0\})$ and $b=\mu(\{+\infty\})$.

This formula can be conveniently rewritten in terms of the
\textit{weighted harmonic mean} $A!_t B$. We define this mean for $0\le t\le 1$ and
positive operators $A$ and $B$, as
\be
A !_t B = (t A^{-1} + (1-t) B^{-1})^{-1}.
\ee
Note that Hansen also defined a weighted harmonic mean, but with a different parametrisation of $t$, ranging over
the interval $[0,+\infty]$ \cite{hansen}.
The extremal cases are $A !_0 B=B$ and $A !_1 B = A$.
In terms of the \textit{parallel sum} \cite{anderson}
\be
A:B := (A^{-1}+B^{-1})^{-1} = B-B(A+B)^{-1}B,
\ee
this formula can be rewritten as
\be\label{eq:whm}
A !_t B = (A/t):(B/(1-t)) = \frac{1}{1-t}\left(B-B\left(\frac{1-t}{t}A+B\right)^{-1}B\right).
\ee
For non-invertible $A$ and/or $B$, one can replace the inverse in the latter formula by the pseudo-inverse.

Performing the substitution $s=(1-t)/t$ (so that $t=1/(1+s)$), we get
$A!_t B = \frac{1+s}{s} \,\,(sA):B$, which appears in the original integral representation (\ref{eq:KAint1}).
Introducing the transformed Radon measure $\nu(t)$ with $d\nu(t) = -d\mu(s)$, we obtain the very simple
representation of a Kubo-Ando mean
\be\label{eq:KAint2}
A\sigma B = \int_{0}^1 \,\,A !_t B\,\,d\nu(t),
\ee
where $aA$ and $bB$ have been absorbed into the integral.
The normalisation condition $1\sigma 1=1$ then imposes the condition $\int_0^1 d\nu(t)=1$, which says that
$d\nu(t)$ must be a probability density. In other words,
the class of
Kubo-Ando means $A\sigma B$ are all possible convex combinations of weighted harmonic means
$A !_t B$, $0\le t\le 1$.

Returning to the axioms that define the Kubo-Ando means, and comparing them to the `dictionary' definition
of means, none of these axioms comes very close in spirit to an in-betweenness property.
The closest match, the monotonicity axiom, is not a comparison between a mean and its arguments,
but a comparison between means of different pairs of arguments.
The conundrum of defining a mean on a partially ordered structure has been solved here in a different way.
Nevertheless, one can still ask the question whether it is not possible to reconcile the two definitions.
Kubo-Ando means might still satisfy an in-betweenness property of some sorts, just like
their classical scalar counterparts, not by definition but as an indirect consequence of its definition.
We will answer this question affirmatively in Section \ref{sec:KA}.
\section{Distance and Angle Monotonicity\label{sec:euclidean}}
In this paper, we shall be dealing with the space of self-adjoint trace class operators.
Endowing this space with the Hilbert-Schmidt inner product $(A,B) = \trace[A^*B]$ turns it
into a real Euclidean vector space.
As positive operators form a subset, $\cS$, of this space,
it makes perfect sense to study $\cS$ from the viewpoint of Euclidean geometry too and consider
Euclidean \textit{distances} and \textit{angles} in $\cS$,
the fact nothwithstanding that nowadays $\cS$ is usually studied from the Riemannian viewpoint, as a manifold
of nonpositive curvature when endowed with the proper metric (see, e.g.\ \cite{bhatia_EMI}).
One of the more obvious benefits of the Euclidean approach is that it also applies to non-positive and
even non-selfadjoint operators.

In accordance with Euclidean geometry, we define Euclidean distance and angles 
in the usual way. These definitions apply, in particular, to positive operators:
\begin{definition}
The Euclidean distance $d$ between two trace class operators $A,B$ is defined as
\be
d(A,B) = \sqrt{\trace[(A-B)^*(A-B)]}.
\ee
\end{definition}

\begin{definition}
The angle $\theta$ between two non-zero trace class operators $A,B$ is defined as
\be
\cos(\theta) = \frac{\Re\trace[A^*B]}{\sqrt{\trace[A^*A]\trace[B^*B]}}.
\ee
\end{definition}
For self-adjoint operators, the $\Re$ operation can obviously be dropped.

Based on the Euclidean distance and angle, we will now define two related in-betweenness properties for
means of positive operators or matrices.
In the following, $A$ and $B$ will always be positive.
We shall say that:
\begin{definition}
An operator mean $\mu$ satisfies \emph{in-betweenness w.r.t.\ Euclidean distance} if and only if 
for all positive $A$ and $B$ the distance between $A$ and $\mu(A,B)$
does not exceed the distance between $A$ and $B$.
\end{definition}
In other words, we shall be demanding that $\mu(A,B)$ lies in the 
Euclidean norm ball with centre $A$ and surface containing $B$.

\begin{definition}
An operator mean $\mu$ satisfies \emph{in-betweenness w.r.t.\ angle} if and only if
for all positive $A$ and $B$ the angle between $A$ and $\mu(A,B)$ does not
exceed the angle between $A$ and $B$.
\end{definition}
This condition requires that $\mu(A,B)$ lies in the 
cone of revolution with origin in the zero operator, 
central axis lying along the direction of $A$, and generated by
the direction of $B$.

One can easily extend these concepts to weighted operator means $\mu(A,B,t)$,
where $t$ is a real scalar in the range $[0,1]$ that expresses how much $A$ dominates over $B$.
The weighted arithmetic mean, for example,
is simply defined as $\mu_{\mbox{ar}}(A,B,t) = tA+(1-t)B$.
In general, $\mu(A,B,0)=B$, $\mu(A,B,1)=A$, and the non-weighted mean is obtained as $\mu(A,B)=\mu(A,B,1/2)$.
For weighted means, the in-betweenness properties can be stated more strongly as monotonicity properties.
\begin{definition}
A weighted operator mean $\mu(A,B,t)$ is \emph{distance-monotonic} if and only if 
the Euclidean distance between $A$ and $\mu(A,B,t)$ decreases monotonically with $t\in[0,1]$.
\end{definition}
Thus, for a distance-monotonic mean, $\trace |A-\mu(A,B,t)|^2$ should decrease monotonically with $t\in[0,1]$.

\begin{definition}
A weighted operator mean $\mu(A,B,t)$ is \emph{angle-monotonic} if and only if 
the angle between $A$ and $\mu(A,B,t)$ decreases monotonically with $t$.
\end{definition}
This condition is equivalent to the monotonic increase of the function
$$
t\mapsto \frac{(\Re\trace A\mu(A,B,t))^2}{\trace\mu(A,B,t)^2}.
$$

In the case of scalar means, distance-monotonicity becomes monotonic decrease of $(a-\mu(a,b,t))^2$,
which is the original in-betweenness property and which should therefore hold for any reasonable 
definition of a weighted scalar mean.
Furthermore, angle-monotonicity is trivially satisfied, as
angles between real positive scalars are always 0.
Finally, it goes without saying that the weighted arithmetic mean is monotonous with respect
to Euclidean distance, because it is geodesic for the Euclidean metric.
\section{Monotonicity of the Power Means\label{sec:power}}
In this section, we will prove that the $p$-power means satisfy in-betweenness, 
both with respect to Euclidean distance
and w.r.t.\ angles, whenever $1\le p\le 2$. Moreover, defining a \textit{weighted $p$-power mean} as
$$
\mu_p(A,B,t)=(t A^p + (1-t)B^p)^{1/p},
$$ 
we show that for $1\le p\le 2$ it is both distance-monotonic and angle-monotonic.

We conjecture that these results holds more generally for larger values of $p$.
The technique we use in our proofs, however, 
ultimately relies on the fact that in the given range of $p$, the function
$x\mapsto x^p$ is convex, while the function $x\mapsto x^{p/2}$ is concave. 
To extend the proofs to larger values of $p$ will require
a different technique.

We begin by showing that for power means 
the monotonicity statement is not really stronger than in-betweenness.
\begin{lemma}
Let $A$ and $B$ be positive operators, and let $f(t)$ be a function of $\mu_p(A,B,t)$ and $A$ (not $B$).
Then $f(t)$ is monotonously increasing over $t\in[0,1]$ if and only if $f(t)\ge f(0)$ for $t\in[0,1]$.
\end{lemma}
\textit{Proof.}
Define $\tilde{B} = \mu_p(A,B,t_1)=(t_1 A^p + (1-t_1)B^p)^{1/p}$ and note that if $t_2\ge t_1$ then
$\mu_p(A,B,t_2)$ can be expressed as a $(t_3 A^p + (1-t_3)\tilde{B}^p)^{1/p}$ for a certain $t_3$ in the interval
$[0,1]$.
Indeed, let $t_3$ be such that $t_2=t_3+(1-t_3)t_1$, then
\beas
\mu_p^p(A,B,t_2) &=& t_2 A^p + (1-t_2)B^p \\
&=&  (t_3+(1-t_3)t_1) A^p +(1-t_3)(1-t_1)B^p \\
&=& t_3 A^p +(1-t_3)\tilde{B}^p = \mu_p^p(A,\tilde{B},t_3).
\eeas
We also have $\mu_p(A,B,t_1) = \tilde{B} = \mu_p(A,\tilde{B},0)$.
Thus the inequality $f(t_1)\le f(t_2)$ reduces to $f(0)\le f(t_3)$ when $B$ is replaced by $\tilde{B}$.
\qed

\begin{theorem}
Let $A$ and $B$ be positive operators,
$0\le t\le 1$ and $1\le p\le2$.
Then $\trace(A-\mu_p(A,B,t))^2$ decreases monotonically with $t$.
\end{theorem}
{\em Proof.}
By the lemma it is enough to show that
\be
\trace(A-\mu_p(A,B,t))^2 \le \trace(A-B)^2. \label{eq:haha}
\ee

Since $\trace(A-B)^2\ge0$,
\be
\trace(A^2+B^2)\ge 2\trace AB. \label{eq:pm1}
\ee
By operator convexity of the function $x\mapsto x^{2/p}$ for $1\le p\le2$,
\be
\trace \mu_p^2(A,B,t) \le \trace \mu_2^2(A,B,t) = \trace(t A^2+(1-t)B^2). \label{eq:pm2}
\ee
Combining (\ref{eq:pm1}), multiplied with $t$, and (\ref{eq:pm2}) gives
\bea
2t\trace AB &\le& \trace(tA^2+tB^2) \nonumber \\
&=& \trace(2tA^2+B^2-(tA^2+(1-t)B^2)) \nonumber \\
&\le& \trace(2tA^2+B^2-\mu_p^2(A,B,t)). \label{eq:pm3}
\eea
By operator concavity of $x\mapsto x^{1/p}$ for $1\le p\le2$,
$$
\mu_p(A,B,t)\ge \mu_1(A,B,t)=tA+(1-t)B,
$$
so that (by the fact that $A\ge0$)
$$
\trace[A \mu_p(A,B,t)] \ge t\trace A^2+(1-t)\trace AB.
$$
Combining this with (\ref{eq:pm3}) gives
\beas
\lefteqn{\trace[\mu_p^2(A,B,t)] - 2\trace[A \mu_p(A,B,t)]} \\
&\le& \trace(2tA^2+B^2-2tAB)-2\trace(tA^2+(1-t)AB) \\
&=& \trace B^2-2\trace AB.
\eeas
Adding $\trace A^2$ to both sides finally gives (\ref{eq:haha}).
\qed

Now we do the same for angle-monotonicity.

\begin{theorem}
Let $A$ and $B$ be positive operators, $0\le t\le 1$ and $1\le p\le2$.
Then the following function of $t$
\be
f(t) := \frac{(\trace[A \mu_p(A,B,t)])^2}{\trace[A^2]\trace[\mu_p(A,B,t)^2]}
\ee
is monotonously increasing with $t$.
\end{theorem}
{\em Proof.}
Again, we can use the lemma to reduce the theorem to the statement $f(0)\le f(t)$
for all $0\le t\le 1$. Dividing out the factor $\trace A^2$ and reorganising the other factors gives:
$$
(\trace[A B])^2\,\trace[\mu_p(A,B,t)^2] \le \trace[B^2]\,(\trace[A \mu_p(A,B,t)])^2.
$$
By absorbing $t$ in $A^p$ and $(1-t)$ in $B^p$,
this is equivalent to
$$
(\trace[A B])^2\,\trace[(A^p+B^p)^{2/p}] \le \trace[B^2]\,(\trace[A (A^p+B^p)^{1/p}])^2.
$$
Let now $a=||A||_2 = (\trace A^2)^{1/2}$ and $b=||B||_2$ and define $G=A/a$ and $H=B/b$. Thus
$G$ and $H$ both have 2-norm equal to 1.
The statement then becomes
$$
(\trace[G H])^2\,\trace[(a^pG^p+b^pH^p)^{2/p}] \le (\trace[G (a^pG^p+b^pH^p)^{1/p}])^2.
$$
Defining $s=a^p/(a^p+b^p)$, which is a convex coefficient, this can be further rewritten as
\be
(\trace[G H])^2\,\trace[(sG^p+(1-s)H^p)^{2/p}] \le (\trace[G (sG^p+(1-s)H^p)^{1/p}])^2.
\label{eq:wa1}
\ee
We will prove this inequality as follows.

First note that the function $x\mapsto x^{2/p}$ is convex, hence
\bea
(\trace[G H])^2\,\trace[(sG^p+(1-s)H^p)^{2/p}] &\le& (\trace[G H])^2\,\trace[sG^2+(1-s)H^2] \nonumber \\
&=& (\trace[G H])^2.\label{eq:wa2}
\eea
Second, the function $x\mapsto x^{1/p}$ is operator concave, hence
$$
sG+(1-s)H \le (sG^p+(1-s)H^p)^{1/p},
$$
so that
\bea
\trace[G (sG^p+(1-s)H^p)^{1/p}] &\ge& \trace[G (sG+(1-s)H)] \nonumber \\
&=& s+(1-s)\trace[G H].\label{eq:wa3}
\eea
Thirdly, by the Cauchy-Schwarz inequality
$$
\trace[G H]\le (\trace[G^2]\trace[H^2])^{1/2} = 1,
$$
so that, for all $0\le s\le 1$,
\be
\trace[G H]\le s+(1-s)\trace[G H]. \label{eq:wa4}
\ee
Combining the three inequalities (\ref{eq:wa2}), (\ref{eq:wa3}) squared, and (\ref{eq:wa4}),
also squared, gives (\ref{eq:wa1}).
\qed
\section{Monotonicity of the Heinz Means\label{sec:heinz}}
In this section we basically prove similar statements as in the previous section but now for the 
(unsymmetrised) Heinz means. As these means are not positive-operator valued, the $\Re$-operation
in the definition of angle is in principle necessary. 
However, it can still be dropped for the 
Heinz means because of their special structure and the fact that $\trace XY$ is real
and positive for positive $X$ and $Y$, even though $XY$ is itself not even Hermitian.

First, we need a simple lemma about convex functions.
\begin{lemma}\label{lem:convex}
Let $x<y$ be real scalars, and
let $a, b$ be distinct real scalars in the open interval $(x,y)$.
When the function $f$ is convex over the interval $[x,y]$, the following holds:
\be
\frac{f(a)-f(x)}{a-x} \le \frac{f(y)-f(b)}{y-b}. \label{eq:xaby}
\ee
\end{lemma}
\textit{Proof.}
Suppose first that $a<b$.
By convexity of $f$ and $a<b<y$, we have $f(b)\le (y-b)f(a)/(y-a) + (b-a)f(y)/(y-a)$,
so that $(f(y)-f(b))/(y-b)\ge (f(b)-f(a))/(b-a)$.
Similarly, from $x<a<b$ follows $(f(b)-f(a))/(b-a)\ge (f(a)-f(x))/(a-x)$.
Combining the two inequalities yields inequality (\ref{eq:xaby}).
For $b<a$ we proceed in a similar way by combining the inequalities
$(f(y)-f(b))/(y-b)\ge (f(a)-f(b))/(a-b)$ and $(f(y)-f(b))/(y-b)\ge (f(a)-f(b))/(a-b)$.
\qed

We start by proving angle-monotonicity for the Heinz means.
\begin{theorem}
Let $A$ and $B$ be positive operators.
For $0\le\nu\le1$,
\be
\trace B^2 (\trace[A^{1+\nu}B^{1-\nu}])^2 \ge \trace[A^{2\nu} B^{2(1-\nu)}](\trace[AB])^2.
\label{eq:heinzth}
\ee
\end{theorem}
\textit{Proof.}
Define the function $x\mapsto g(x) = \trace[A^{2x} B^{2(1-x)}]$. Inequality (\ref{eq:heinzth}) says that $g(x)$
obeys
$$
g(0) g^2(1/2+\nu/2) \ge g(\nu)g^2(1/2).
$$
Upon taking logarithms and rearranging terms, this is equivalent to
$$
\log g(\nu)-\log g(0) \le 2(\log g(1/2+\nu/2)-\log g(1/2)).
$$
This would follow from Lemma \ref{lem:convex}, with $x=0$, $y=1/2+\nu/2$, $a=1/2$ and $b=\nu$, provided $\log g(x)$
is convex.
This convexity is now easily seen to be equivalent with a Cauchy-Schwarz inequality:
$$
(\trace[A^{x+y} B^{2-(x+y)}])^2 \le \trace[A^{2x} B^{2(1-x)}]\,\trace[A^{2y} B^{2(1-y)}].
$$
Taking logarithms gives the required statement
$$
\log g((x+y)/2) \le (\log g(x)+\log g(y))/2.
$$
\qed

The corresponding result for distance-monotonicity is proven in a similar way.
\begin{theorem}
Let $A$ and $B$ be positive semidefinite matrices.
For $0\le\nu\le1$,
\be
\trace B^2 +2\trace[A^{1+\nu}B^{1-\nu}] \ge \trace[A^{2\nu} B^{2(1-\nu)}]+2\trace[AB].
\label{eq:heinzth2}
\ee
\end{theorem}
\textit{Proof.}
The proof proceeds in the same way as before, but now exploiting the convexity of $g(x)$.
The latter follows immediately
from the convexity of $\log g(x)$ by the fact that $\exp(x)$ is a monotonously increasing convex function.
\qed
\section{Monotonicity of the Kubo-Ando Means\label{sec:KA}}
Given the initial success in finding two non-trivial operator means for which in-betweenness holds,
it would be very interesting if the larger class of Kubo-Ando means also satisfied it.
This, however, is not the case; at least, not with the current definition of in-betweenness.

We will consider a simple counterexample for the harmonic mean. As already stated,
the Kubo-Ando means are built up from the harmonic mean.
We choose the following $2\times 2$ matrices:
$$
A=\twomat{5}{7}{7}{10},\qquad
B=\twomat{5}{2}{2}{1}.
$$
A simple numerical calculation reveals that,
for $t$ between $0$ and about $0.32$, the Euclidean distance between $A$ and $A!_t B$ increases, rather than
decreases with $t$. Thus for $t$ in the interval $[0,0.32]$, $||A-A!_tB||_2 \not\le ||A-B||_2$.
By rescaling $A$ and $B$ one can make this happen at any value of $t$, including $t=1/2$.

This appears to be very unsatisfying at first, but then one has to realise that this is all a matter of geometry.
One cannot really expect that quantities that are natural in one geometry should possess properties belonging to another.
The strong connections between Kubo-Ando means and hyperbolic Riemannian geometry (non-positive curvature) suggest
that maybe one should modify the definition of in-betweenness to reflect this different geometry, 
and more particularly the distance measure used to define in-betweenness.
This is what we will attempt in the following paragraphs. 

Given that the Kubo-Ando means are convex combinations of weighted harmonic means $!_t$, 
one might start
with defining an in-betweenness that works for all the weighted harmonic means, and then take it from there.
Looking at the defining formula (\ref{eq:whm}), the following candidate for a distance measure comes to mind:
$d_{-1}(A,B) = ||A^{-1}-B^{-1}||_2$, which is the Euclidean distance between the inverses.
It is trivial to see that, with this candidate, in-betweenness holds for all weighted harmonic means.
Indeed, for invertible, positive $A$ and $B$, and $0\le t\le 1$:
\beas
d_{-1}(A,A!_t B) &=& ||A^{-1}-(A!_t B)^{-1} ||_2  \\
&=& ||A^{-1}-(tA^{-1}+(1-t) B^{-1}) ||_2 \\
&=& (1-t)\,\,||A^{-1}-B^{-1}||_2 = (1-t)\,\, d_{-1}(A,B).
\eeas
This actually shows that $t\mapsto A!_tB$ is a minimal geodesic with respect to the 
inverted Euclidean metric $d_{-1}(A,B)$, just like the weighted arithmetic mean defines a geodesic
w.r.t.\ the ordinary Euclidean metric $d$.

Now does it also work for general convex combinations of the
weighted harmonic means? The answer is no -- it cannot, for the following reason. 
Given any Kubo-Ando mean $\sigma$, we get another Kubo-Ando mean $\sigma^*$
called the adjoint via the correspondence $A\sigma^* B = (A^{-1} \sigma B^{-1})^{-1}$. Hence
a Kubo-Ando mean satisfies in-betweenness w.r.t.\ inverted Euclidean distance 
$d_{-1}$ if
and only if its adjoint mean satisfies in-betweenness w.r.t.\ the ordinary Euclidean distance $d$. 
The latter condition is not always satisfied, as shown by the counterexample above.

This suggests that to find a distance measure for which all Kubo-Ando means are monotonic one should
be looking at a distance for which the geodesic is in the `middle' of the set of means.
One obvious candidate is the set
of weighted geometric means $A \#_t B = A^{1/2} (A^{-1/2} B A^{-1/2})^t A^{1/2}$.
It is well-known that these means
define a geodesic $t\mapsto A\#_t B$, with respect to the so-called \textit{trace metric distance}
$\delta(A,B) = ||\log(A^{-1/2}BA^{-1/2})||_2$. 
That is, $\delta(A,A\#_t B) = t \delta(A,B)$.
Hence these means trivially satisfy in-betweenness w.r.t.\ the trace metric distance.

It turns out that all Kubo-Ando means are monotonic w.r.t.\ $\delta$.
\begin{theorem}
For all positive operators $A,B$, and Kubo-Ando means $\sigma$,
\be
\delta(A,A\sigma B) \le \delta(A,B).
\ee
\end{theorem}
\textit{Proof.}
We exploit the fact that $\delta$ is invariant under conjugations. That is, for all $M$,
$\delta(MAM^*,MBM^*)=\delta(A,B)$. Writing $B=A^{1/2}CA^{1/2}$, we then see
that we only have to prove the statement for $A=\id$ and $B=C$.

Let $c$ be a real scalar, with $0\le c\le 1$.
Then we have, for $0\le t\le 1$,
$$
1\le t+(1-t)c^{-1}\le c^{-1}.
$$
Inverting and taking the integral over $t$ with probability measure $dp(t)$ gives
$$
1\ge \int_0^1 dp(t) (t+(1-t)c^{-1})^{-1}\ge c.
$$
As every Kubo-Ando mean $\sigma$ can be written as a convex combination of weighted harmonic means, this shows that
for every such mean, $1\ge 1\sigma c \ge c$.

Since the function $|\log x|$ is monotonically decreasing for $x\le 1$, this implies
$$
|\log 1\sigma c| \le |\log c|.
$$
The same inequality can be shown to hold when $0< c\le 1$.
Passing to positive operators $C$ and taking the $||.||_2$ norm gives the required inequality
$$
\delta(\id,\id\sigma C)=||\log(\id\sigma C)||_2 \le ||\log C||_2=\delta(\id,C).
$$
\qed

From the proof one sees that this theorem holds 
more generally for every mean $\mu$ that satisfies the closure axiom (to have positivity), achieves
equality in the transformer inequality (to be able to apply the invariance of $\delta$ under conjugations) and
scalar in-betweenness $x\le \mu(x,y)\le y$.

\section{Conclusion\label{sec:discussion}}
In this paper we have introduced the notion of in-betweenness for operator means, and the slightly stronger
one of monotonicity with respect to a given distance measure (or metric), for those operator means that admit 
a simple parameterisation. These notions can be seen as a relaxation of geodesity, in the following sense. 
When a parameterised
operator mean $\mu(x,y,t)$ traces out a geodesic $t\mapsto \mu(x,y,t)$ with respect to a given metric $d$,
it satisfies $d(\mu(x,y,t),y)=(1-t)d(x,y)$ by definition. 
This, however, requires a careful matching between the parameterisation
of the mean and the chosen metric. This may not always be possible, be it for internal or for external reasons. 
In that case it might still be useful to have monotonicity, 
which is the \textit{inequality} $d(\mu(x,y,t),y)\le (1-t)d(x,y)$.
We have exhibited two non-trivial examples of operator means that are monotonous with respect to
the Euclidean metric, even though the Euclidean metric would not be considered the natural one for 
these means. We have also shown that all the Kubo-Ando means are monotonous w.r.t.\ the trace metric; in contrast
only the (weighted) geometric means are geodesic in this metric.

In this work we have only scratched the surface and many questions remain. 
Most importantly, it would be very interesting if one could give
a full characterisation of operator means that are monotonic w.r.t.\ a given metric, and possibly come up
with an alternative axiomatic approach to operator means.
\section{Acknowledgments}
I am grateful for the hospitality of the Institut Mittag-Leffler, Djursholm (Sweden) where this manuscript
was completed.

\end{document}